# THE SEMI-COMMUTATOR OF TOEPLITZ OPERATORS ON THE BIDISC

CAIXING GU AND DECHAO ZHENG

ABSTRACT. In this paper we characterize when the semi-commutator $T_f T_g - T_{fg}$ of two Toeplitz operators $T_f$ and $T_g$ on the Hardy space of the bidisc is zero. We also show that there is no nonzero finite rank semi-commutator on the bidisc. Furthermore explicit examples of compact semi-commutators with symbols continuous on the bitorus $T^2$ are given.

## 1. INTRODUCTION

Let D be the open unit disk in $C$. Its boundary is the unit circle $T$. The bidisc $D^2$ and the torus $T^2$ are the subsets of $C^2$ which are Cartesian products of two copies $D$ and $T$, respectively. Let $d\sigma(z)$ be the normalized Haar measure on $T^2$. The Hardy space $H^2(D^2)$ is the closure of the polynomials in $L^2(T^2, d\sigma)$ (or $L^2(T^2)$). Let $P$ be the orthogonal projection from $L^2(T^2)$ onto $H^2(D^2)$. The Toeplitz operator with symbol $f$ in $L^\infty(T^2)$ is defined by $T_f(h) = P(fh)$, for all $h \in H^2(D^2)$ and the Hankel operator with symbol $f$ is defined by $H_f(h) = (I - P)(fh)$, for all $h \in H^2(D^2)$.

Let $f$ and $g$ be two bounded functions on $T^2$. In this paper we study the the semi-commutator $T_f T_g - T_{fg}$ of two Toeplitz operators $T_f$ and $T_g$ on the bidisc. As in the case of the unit disk, the semi-commutator is connected to the Hankel operators by the following relation.

$$(1) \qquad T_f T_g - T_{fg} = -H_{\bar{f}}^* H_g.$$

To motivate the problems to be considered in the paper, we shall recall some classical results for the semi-commutators of Toeplitz operators on the Hardy space $H^2(D)$ of the unit disk. Let $f_1$ and $g_1$ be bounded functions on the unit circle $T$. For two Toeplitz operators $T_{f_1}$ and $T_{g_1}$ on the Hardy space $H^2(D)$ of the unit disk, Brown-Halmos [2] shows that the semi-commutator $T_{f_1} T_{g_1} - T_{f_1 g_1}$ is zero if and only if either $\bar{f}_1$ or $g_1$ is analytic. In other words, $H_{\bar{f}_1}^* H_{g_1}$ is zero if and only if either $H_{\bar{f}_1}$ or $H_{g_1}$ is zero.

In this paper, we will characterize when the semi-commutator $T_f T_g - T_{fg}$ of two Toeplitz operators $T_f$ and $T_g$ on the bidisc is zero. In particular, we note that unlike

1991 *Mathematics Subject Classification.* Primary 47B35.

Gu's research is supported in part by NSF Grant DMS-9022140 during residence at MSRI. Zheng is supported in part by a NSF grant.

on the unit disk one can have both $H_{\bar{f}}$ and $H_g$ are not zero, but their product $H_{\bar{f}}^* H_g$ is zero. Furthermore we will see that there is no finite rank semi-commutator on the bidisc. But this is false on the unit disk. Indeed it was shown in [1] that for Toeplitz operators $T_{f_1}$ and $T_{g_1}$ on $H^2(D)$, $T_{f_1} T_{g_1} - T_{f_1 g_1}$ is of finite rank if and only if either $\bar{f}_1$ or $g_1$ is an analytic function plus a rational function.

The main question to be considered here is when the semi-commutator $T_f T_g - T_{fg}$ on the bidisc is compact. This problem is connected with the spectral theory of Toeplitz operators on the bidisc and various applications, see [7], [8] and reference therein.

For Toeplitz operators $T_{f_1}$ and $T_{g_1}$ on the unit disk, this hard problem was solved by the combined efforts of Axler, Chang, and Sarason [1] and Volberg [12]. Their beautiful result is that that $T_{f_1 g_1} - T_{f_1} T_{g_1}$ is compact if and only if $H^\infty[\bar{f}_1] \bigcap H^\infty[g_1] \subset H^\infty(D) + C(T)$; here $H^\infty[g_1]$ denotes the closed subalgebra of $L^\infty(T)$ generated by $H^\infty(D)$ and $g_1$ and $C(T)$ the continuous functions on $T$.

The function theory on the bidisc is quite different from and much less understood than the function theory on the unit disk [9], [3] and [5]. The proof of the above result on the unit disk relies on some deep results and techniques from function theory on the unit disk which are not available from function theory on the bidisc.

In this paper we content ourselves with some partial results for the compactness of $T_f T_g - T_{fg}$. By carefully analyzing the action of $T_f T_g - T_{fg}$ on the reproducing kernel of $H^2(D^2)$ and exploiting the harmonicity of certain functions, we will get a necessary condition for the compactness of $T_f T_g - T_{fg}$. This shows that for a large class of functions $f$ and $g$, $T_f T_g - T_{fg}$ is compact if and only if it is zero. For example, if $f$ is a trigonometric polynomial on $T^2$ and $g$ is an arbitrary bounded function on $T^2$, then for such $f$ and $g$, there is no compact semi-commutator $T_f T_g - T_{fg}$. Also as a corollary of this condition we see that there are no compact Hankel operators on bidisc, which was proved by M. Cotlar and C. Sadosky [4] using a completely different method. It is natural to guess that for $f$ and $g$ in $L^\infty(T^2)$, $T_f T_g - T_{fg}$ is compact if and only if $T_f T_g - T_{fg}$ is zero. But it is false.

For the class of bounded functions f and g on $T^2$ of the form $f(z_1, z_2) = f_1(z_1) f_2(z_2)$ and $g(z_1, z_2) = g_1(z_1) g_2(z_2)$, we will show that $T_f T_g - T_{fg}$ is compact if and only if the following two conditions hold.

(1) $f_1(z_1) g_1(z_1) = 0$ and $f_2(z_2) g_2(z_2) = 0$ on $T$, and

(2) $\lim_{z \to T^2} \|H_{\overline{f}} k_z\|_2 \|H_g k_z\|_2 = 0$,

where $k_z$ denotes the normalized reproducing kernel of $H^2(D^2)$. Here the Littlewood-Paley theory on the bidisc and a certain distribution function inequality in Zheng [13] play a key role. We remark that only one condition similar to Condition (2) above is needed for the compactness of the semi-commutator $T_{f_1} T_{g_1} - T_{f_1 g_1}$ on the Hardy space $H^2(D)$ of the unit disk. Indeed it was shown in [13] that $T_{f_1} T_{g_1} - T_{f_1 g_1}$ is compact on $H^2(D)$ if and only if $\lim_{z_1 \to T} \|H_{\overline{f_1}} k_{z_1}\|_2 \|H_{g_1} k_{z_1}\|_2 = 0$, where $k_{z_1}$ denotes the normalized reproducing kernel of $H^2(D)$.

Now we outline the plan of the paper. In Section 2 we study when the semi-commutator of two Toeplitz operators is zero or finite rank. In Section 3 we derive a necessary condition for the compactness of $T_f T_g - T_{fg}$. In Section 4, for f and g of the form $f(z_1, z_2) = f_1(z_1)f_2(z_2)$ and $g(z_1, z_2) = g_1(z_1)g_2(z_2)$, we characterize when $T_f T_g - T_{fg}$ is compact. This result shows that the necessary condition for the compactness of $T_f T_g - T_{fg}$ obtained in Section 3 is not sufficient in general. This result also provides us with explicit examples of compact semi-commutators $T_f T_g - T_{fg}$ with symbols $f$ and $g$ continuous on the bitorus $T^2$.

## 2. Zero semi-commutators

Let $Z$ denote the set of all integers, $Z_+$ the set of all nonnegative integers and $Z_-$ the set of all negative integers. As in [10] we consider multiple Fourier series on the bitorus $T^2$. The multiple Fourier series on the bitorus $T^2$ can be viewed as the Fourier transformation on $L^1(T^2)$. For f in $L^1(T^2)$, the Fourier transformation of $f$ on $Z \times Z$ is given by

$$f_m = (\frac{1}{2\pi})^2 \int_0^{2\pi} \int_0^{2\pi} f(e^{i\theta_1}, e^{i\theta_1}) e^{i(m,\theta)} d\theta_1 d\theta_2$$

where $m = (m_1, m_2) \in Z \times Z$, $\theta = (\theta_1, \theta_2)$ and $(m, \theta) = m_1\theta_1 + m_2\theta_2$. By Theorem 1.7 in [SW], the Fourier transformation is injective, i.e. If $f \in L^1(T^2)$ and $f_m = 0$ for all $m \in Z \times Z$, then $f \equiv 0$. We recall also that by using multiple Fourier series,

$$L^2(T^2) = \{f : f = \sum_{m \in Z \times Z} f_m e^{i(m,x)}, \sum_{m \in Z \times Z} |f_m|^2 < \infty\}$$

$$H^2(D^2) = \{h : h = \sum_{m \in Z_+ \times Z_+} h_m e^{i(m,x)}, \sum_{m \in Z_+ \times Z_+} |h_m|^2 < \infty\}$$

and

$$Pf = \sum_{m \in Z_+ \times Z_+} f_m e^{i(m,x)} \text{ for } f = \sum_{m \in Z \times Z} f_m e^{i(m,x)} \in L^2(T^2).$$

**Theorem 1.** *Let $f$ and $g$ be two bounded functions on the torus $T^2$. The following are equivalent.*
*(1) $T_f T_g - T_{fg}$ is a finite rank operator.*
*(2) $T_f T_g - T_{fg}$ is zero.*
*(3) For each $i$ $(i = 1, 2)$ either $\bar{f}$ or $g$ is analytic in $z_i$.*

**Proof.** We first show that (3) implies (2). Assume that for each $i$ $(i = 1, 2)$ either $\bar{f}$ or $g$ is analytic in $z_i$. Without loss of generality, assume that $\bar{f}$ is analytic in $z_1$ for and $g$ is analytic in $z_2$. Then it is easy to see that

$$(I - P)(gh_1) = \sum_{m=(m_1,m_2) \in Z_- \times Z_+} a_m z^m$$

for all $h_1 \in H^2(D^2)$ and

$$(I - P)(\bar{f}h_2) = \sum_{m=(m_1, m_2) \in Z_+ \times Z_-} b_m z^m$$

for all $h_2 \in H^2(D^2)$. Therefore

$$(H_g h_1, H_{\bar{f}} h_2) = 0$$

for all $h_1, h_2 \in H^2(D^2)$. That is, $H_{\bar{f}}^* H_g = 0$. So $T_f T_g - T_{fg} = -H_{\bar{f}}^* H_g = 0$

It is obvious that (2) implies (1). Now we prove that (1) implies (3). That is, assume that $T_f T_g - T_{fg}$ is a finite rank operator, we will show that for each $i$ ($i = 1, 2$) either $\bar{f}$ or $g$ is analytic in $z_i$. Without loss of generality, we assume that $i = 1$. We write $f$ and $g$ as

$$\bar{f} = \sum_{i=-\infty}^{\infty} f_i(z_2) z_1^i = \sum_{i=-\infty}^{\infty} f_i z_1^i$$

$$g = \sum_{i=-\infty}^{\infty} g_i(z_2) z_1^i = \sum_{i=-\infty}^{\infty} g_i z_1^i$$

Let $\alpha, \beta, k, l \in Z_+$. Then by a straightforward computation we have

$$
\begin{aligned}
& (H_f^* H_g z_1^k z_2^\alpha, z_1^l z_2^\beta) \\
= & (H_g z_2^\alpha z_1^k, H_f z_2^\beta z_1^l) \\
= & \sum_{i \leq -(k+1)} (g_i(z_2) z_2^\alpha, f_{-(i+k-l)}(z_2) z_2^\beta) + \sum_{i \geq -k} ((I_2 - P_2) g_i(z_2) z_2^\alpha, f_{-(i+k-l)}(z_2) z_2^\beta)
\end{aligned}
$$

where $I_2$ is the identity on $L^2(T)$ and $P_2$ is the projection from $L^2(T)$ onto $H^2(D)$. Therefore

$$
\begin{aligned}
& (H_g z_1^k z_2^\alpha, H_f z_1^l z_2^\beta) - (H_g z_1^{k+1} z_2^\alpha, H_f z_1^{l+1} z_2^\beta)) \\
= & (g_{-(k+1)}(z_2)(z_2)^\alpha, f_{-(l+1)}(z_2)(z_2)^\beta) - ((I_2 - P_2) g_{-(k+1)}(z_2)(z_2)^\alpha, f_{-(l+1)}(z_2)(z_2)^\beta) \\
= & (P_2 g_{-(k+1)}(z_2)(z_2)^\alpha, f_{-(l+1)}(z_2)(z_2)^\beta).
\end{aligned}
$$

Let $S_1$ denote the multiplication by $z_1$ on $H^2(D^2)$, i.e., $S_1 h = z_1 h$ for $h \in H^2(D^2)$. The above relation implies that

$$(S_1^{*l} H_f^* H_g S_1^k - S_1^{*(l+1)} H_f^* H_g S_1^{k+1}) h_2(z_2) = T_{\overline{f_{-(l+1)}}} T_{g_{-(k+1)}} h_2(z_2)$$

for all $h_2 \in H^2(D)$. Therefore, if $H_f^* H_g$ is a finite rank operator on $H^2(D^2)$, then $T_{\overline{f_{-(l+1)}}} T_{g_{-(k+1)}}$ is a finite rank operator on $H^2(D)$. By a result in [2], we have that either $f_{-(l+1)}$ or $g_{-(k+1)} = 0$. Hence either $f_{-(l+1)} = 0$ for all $k \geq 0$ or $g_{-(k+1)} = 0$ for all $l \geq 0$. That is either $\bar{f}$ or $g$ is analytic in $z_1$. This finishes the proof of the theorem.

### 3. A NECESSARY CONDITION FOR COMPACTNESS

In this section we will give a necessary condition for the compactness of the semi-commutator $T_f T_g - T_{fg}$.

Before going to the main result of this section, we need some notations and definitions. We use $z$ to denote the vector $(z_1, z_2)$ in $C^2$ of two dimensional complex plane. Since for any z in $D^2$, the pointwise evaluation of functions in $H^2(D^2)$ at z is a bounded functional, there is a function $K_z$ in $H^2(D^2)$ such that

$$f(z) = (f, K_z)$$

for all f in $H^2(D^2)$. $K_z(w)$ is called the reproducing kernel for $H^2(D^2)$ and sometimes we use K(z,w) to denote $K_z(w)$.

Let $K_{z_1}$ denote the reproducing kernel of the Hardy space $H^2(D)$ of the unit disk at the point $z_1$ in $D$, and $k_{z_1}$ the normalized reproducing kernel of the Hardy space $H^2(D)$ at the point $z_1 \in D$. Namely,

$$K_{z_1} = \frac{1}{(1 - \overline{z_1} w_1)}, \quad k_{z_1} = \frac{(1 - |z_1|^2)^{1/2}}{(1 - \overline{z_1} w_1)}.$$

It is easy to check that the reproducing kernel $K_z$ of the Hardy space $H^2(D^2)$ of the bidisc is the product $K_{z_1}(w_1) K_{z_2}(w_2)$ of the reproducing kernels of $H^2(D)$. So the normalized reproducing kernel $k_z$ of $H^2(D^2)$ is also the product $k_{z_1}(w_1) k_{z_2}(w_2)$ of the normalized reproducing kernels of $H^2(D)$. We observe that $k_z$ weakly converges to zero in $H^2(D^2)$ as $z$ tends to the boundary of $D^2$. For a function $f$ in the space $L^2(T^2)$, i.e.,

$$f(w) = \sum_{m \in Z \times Z} f_m w^m, \quad w \in T^2,$$

where $f_m$ is a sequence of numbers such that

$$\sum_{m \in Z \times Z} |f_m|^2 < \infty,$$

the bi-harmonic (in short, harmonic) extension $f(z)$ of $f$ is defined via the Possion integral

$$f(z) := \int_T \int_T f(w) \frac{(1 - |z_1|^2)}{|1 - \overline{z_1} w_1|^2} \frac{(1 - |z_2|^2)}{|1 - \overline{z_2} w_2|^2} d\sigma(w_1) d\sigma(w_2)$$

$$= (f(w) k_z(w), k_z(w)).$$

We write the power series expansion of the harmonic extension $f(z)$ of $f$ as follows:

$$f(z) = \sum_{m \in Z \times Z} f_m z^m$$

$$\text{(2)} \qquad := f_{++}(z) + f_{+-}(z) + f_{-+}(z) + f_{--}(z),$$

where

$$f_{++}(z) := \sum_{m \in Z_+ \times Z_+} f_m z^m, \quad f_{+-}(z) := \sum_{m \in Z_+ \times Z_-} f_m z^m$$

$$f_{-+}(z) := \sum_{m \in Z_- \times Z_+} f_m z^m, \quad f_{--}(z) := \sum_{m \in Z_- \times Z_-} f_m z^m.$$

Also let

$$(3) \qquad f_+(z) = f_{++}(z) + f_{+-}(z), \quad f_-(z) = f_{-+}(z) + f_{--}(z),$$

where for example, $m = (m_1, m_2) \in Z_- \times Z_+$ means that $m_1 \in Z_+$ and $m_2 \in Z_-$, $z^m$ the product $\bar{z}_1^{|m_1|} z_2^{m_2}$.

**Theorem 2.** *Let $f$ and $g$ be two bounded functions on the bitorus $T^2$. Then $T_f T_g - T_{fg}$ is compact implies that the following two statements hold.*
*(1) For all $z_1 \in D$ and all $z_2 \in T$,*

$$\frac{\partial f}{\partial z_1}(z_1, z_2) \frac{\partial g}{\partial \bar{z}_1}(z_1, z_2) = 0.$$

*(1) For all $z_2 \in D$ and all $z_1 \in T$,*

$$\frac{\partial f}{\partial z_2}(z_1, z_2) \frac{\partial g}{\partial \bar{z}_2}(z_1, z_2) = 0.$$

**Proof.** Without loss of generality, we prove statement (1). The proof of statement (2) is similar. We write $f$ and $g$ as (see above for the notations)

$$f(z) := f_{++}(z) + f_{+-}(z) + f_{-+}(z) + f_{--}(z)$$

$$:= f_+(z) + f_-(z)$$

and similarly

$$g(z) := g_{++}(z) + g_{+-}(z) + g_{-+}(z) + g_{--}(z)$$

$$:= g_+(z) + g_-(z).$$

We compute the action of a Toeplitz operator $T_g$ for $g$ as above on the normalized reproducing kernel $k_z(w)$ as follows:

$$(4) \qquad \begin{aligned} & T_g k_z(w_1, w_2) \\ = & g_{++}(w_1, w_2) k_z + g_{+-}(w_1, z_2) k_z + g_{-+}(z_1, w_2) k_z + g_{--}(z_1, z_2) k_z \end{aligned}$$

for $z = (z_1, z_2)$ in $D^2$, where for example

$$f_{+-}(w_1, z_2) = \sum_{(m_1, m_2) \in Z_+ \times Z_-} g_m w_1^{m_1} \bar{z}_2^{|m_2|}$$

for all $w_1 \in T$ and $z_2 \in D$. We remark that in fact the above formula (4) holds for any $g \in L^2(T^2)$.

For convenience, we consider $T_{\bar{f}}T_g - T_{\bar{f}g}$ instead of $T_f T_g - T_{fg}$. We write that

$$((T_{\bar{f}}T_g - T_{\bar{f}g})k_z, k_z)$$
$$= -\int_{T^2} \overline{f(w)}g(w)|k_z(w)|^2 d\sigma(w) + (T_g k_z, T_f k_z)$$
$$(5) \qquad := h_1(z) + h(z)$$

for any z in $D^2$, where

$$h_1(z) = -\int_{T^2} \overline{f(w)}g(w)|k_z(w)|^2 d\sigma(w),$$

$$h(z) = (T_{\bar{f}}T_g k_z, k_z).$$

Furthermore by (5), we write $h(z)$ as

$$h(z) = (T_{\bar{f}}T_g k_z, k_z) = (T_{\bar{f}}T_g k_z, k_z)$$
$$:= h_{11}(z) + h_{12}(z) + h_{21}(z) + h_{22}(z)$$

where

$$(6) \qquad h_{11}(z) := (T_{g_+} k_z(w), T_{f_+} k_z(w)), \quad h_{12}(z) := (T_{g_+} k_z(w), T_{f_-} k_z(w)),$$

$$(7) \qquad h_{21}(z) := (T_{g_-} k_z(w), T_{f_+} k_z(w)), \quad h_{22}(z) := (T_{g_-} k_z(w), T_{f_-} k_z(w)).$$

Let $m = (m_1, m_2) \in Z \times Z$ be fixed. Let $\theta = (\theta_1, \theta_2)$ and

$$U_\theta = diag\{e^{i\theta_1}, e^{i\theta_2}\}.$$

Set

$$H_1(z) := \int_{T^2} h_1(U_\theta z) e^{i(m,\theta)} d\theta$$

$$H_{ij}(z) := \int_{T^2} h_{ij}(U_\theta z) e^{i(m,\theta)} d\theta, \quad i,j = 1,2$$

for $z \in D^2$, where $h_1(U_\theta z)$ and $h_{ij}(U_\theta z)$ are obtained by replacing $z$ with $U_\theta z$ in the definition of $h_1$ and $h_{ij}$ respectively. For example, by (4), we have

$$H_{12}(z) := \int_{T^2} (g_{++}(w_1, w_2)k_{U_\theta z} + g_{+-}(w_1, e^{i\theta_2}z_2)k_{U_\theta z},$$
$$(8) \qquad f_{-+}(e^{i\theta_1}z_1, w_2)k_{U_\theta z} + f_{--}(e^{i\theta_1}z_1, e^{i\theta_2}z_2)k_{U_\theta z})e^{i(m,\theta)} d\theta.$$

For a fixed $z_1 \in D$ and $u_2 \in T$, let $z_{2\alpha}$ converge to the point $u_2$ on $T$, then $z_\alpha = (z_1, z_{2\alpha})$ converges to the boundary point $(z_1, u_2)$ of the bidisc $D^2$. We claim that $H_1(z)$, $H_{11}(z)$, $H_{12}(z)$, $H_{21}(z)$ and $H_{22}(z)$ are continuous on $\overline{D^2}$; and furthermore if we denote the limits of $H_1(z)$, $H_{11}(z)$, $H_{12}(z)$, $H_{21}(z)$ and $H_{22}(z)$ as $z_\alpha$ converges to the boundary point $(z_1, u_2)$ by $H_1(z_1, u_2)$, $H_{11}(z_1, u_2)$, $H_{12}(z_1, u_2)$, $H_{21}(z_1, u_2)$ and $H_{22}(z_1, u_2)$, respectively. Then
$H_1(z_1, u_2)$, $H_{11}(z_1, u_2)$, $H_{12}(z_1, u_2)$ and $H_{21}(z_1, u_2)$ are harmonic in $z_1$ for any $u_2 \in T$. However $H_{22}(z_1, u_2)$ is not necessary harmonic in $z_1$. We postpone the proof of the claim and first continue with the proof of the theorem.

We will first show that if $T_{\bar{f}}T_g - T_{\bar{f}g}$ is compact, then $H_{22}(z)$ is harmonic in $z_1$. Replacing z by $U_\theta z$ in (5) yields

$$((T_{\bar{f}}T_g - T_{\bar{f}g})k_{U_\theta z}, k_{U_\theta z}) = h_1(U_\theta z) + h(U_\theta z).$$

Multiplying the above equation by $e^{i(m,\theta)}$ and then integrating with respect to $\theta$ give that

$$\int_{T^2} ((T_{\bar{f}}T_g - T_{\bar{f}g})k_{U_\theta z}, k_{U_\theta z})e^{i(m,\theta)}d\theta$$

$$= H_1(z) + H_{11}(z) + H_{12}(z) + H_{21}(z) + H_{22}(z).$$

Now note that $k_z$ weakly converges to zero as z goes to the boundary of $D^2$. Therefore if $T_{\bar{f}}T_g - T_{\bar{f}g}$ is compact, then

$$\lim_{(z)_\alpha \to (z_1, u_2)} \int_{T^2} ((T_{\bar{f}}T_g - T_{\bar{f}g})k_{U_\theta z_\alpha}, k_{U_\theta z_\alpha})e^{i(m,\theta)}d\theta = 0.$$

Hence by (5), we conclude that

$$\lim_{z_\alpha \to (z_1, u_2)} (H_1(z) + H_{11}(z) + H_{12}(z) + H_{21}(z) + H_{22}(z)) = 0$$

for any $z_1 \in D$ and $u_2 \in T$. That is

$$-H_{22}(z_1, u_2) = H_1(z_1, u_2) + H_{11}(z_1, u_2) + H_{12}(z_1, u_2) + H_{21}(z_1, u_2).$$

Since by our claim $H_1(z_1, u_2)$, $H_{11}(z_1, u_2)$, $H_{12}(z_1, u_2)$ and $H_{21}(z_1, u_2)$ are harmonic in $z_1$ for any $u_2 \in T$, so $H_{22}(z_1, u_2)$ is harmonic in $z_1$ for any $u_2 \in T$. Thus $\Delta_{z_1} H_{22}(z_1, u_2) = 0$.

On the other hand, if we write

$$f_-(z_1, z_2) = \sum_{k=1}^\infty \frac{1}{k!} \frac{\partial^k f_-}{\partial \bar{z}_1^k}(0, z_2)\bar{z}_1^k$$

and

$$g_-(z_1, z_2) = \sum_{k=0}^\infty \frac{1}{k!} \frac{\partial^k g_-}{\partial \bar{z}_1^k}(0, z_2)\bar{z}_1^k.$$

Let $m = (0, m_2)$. We note that

$$H_{22}(0, u_2) =$$

$$\sum_{k=1}^\infty \frac{1}{(k!)^2} \int_T \overline{\frac{\partial^k f_-}{\partial \bar{z}_1^k}(0, u_2 e^{i\theta_2})} \frac{\partial^k g_-}{\partial \bar{z}_1^k}(0, u_2 e^{i\theta_2})e^{im_2\theta_2}d\theta_2 |z_1|^{2k}.$$

Since $H_{22}(z_1, u_2)$ is harmonic in $z_1$ for any $u_2 \in T$. Thus $\Delta_{z_1} H_{22}(z_1, u_2) = 0$. In particular, $\Delta_{z_1} H_{22}(z_1, u_2) = 0$ for $z_1 = 0$. Hence

$$\int_T \overline{\frac{\partial f_-}{\partial \bar{z}_1}(0, u_2 e^{i\theta_2})} \frac{\partial g_-}{\partial \bar{z}_1}(0, u_2 e^{i\theta_2}) e^{im_2\theta_2} d\theta_2 = 0$$

for all $m_2 \in Z$ and $u_2 \in T$. We observe that both $\dfrac{\partial f_-}{\partial \bar{z}_1}(0, u_2)$ and $\dfrac{\partial g_-}{\partial \bar{z}_1}(0, u_2)$ are in $L^2(T)$. So $\overline{\dfrac{\partial f_-}{\partial \bar{z}_1}(0, u_2)} \dfrac{\partial g_-}{\partial \bar{z}_1}(0, u_2)$ is in $L^1(T)$ and the Fourier transformation of $\overline{\dfrac{\partial f_-}{\partial \bar{z}_1}(0, u_2)} \dfrac{\partial g_-}{\partial \bar{z}_1}(0, u_2)$ on $Z$ is zero. The injection of the Fourier transformation implies that

$$\overline{\frac{\partial f_-}{\partial \bar{z}_1}(0, u_2)} \frac{\partial g_-}{\partial \bar{z}_1}(0, u_2) = 0$$

for all $u_2 \in T$. But this is same as

$$\overline{\frac{\partial f}{\partial \bar{z}_1}(0, u_2)} \frac{\partial g}{\partial \bar{z}_1}(0, u_2) = 0$$

for all $u_2 \in T$.

Next by using the Möbius transform of the bidisc, we will show that statement (1) in the theorem holds. Let $\phi_z(w)$ denote the Möbius transform

$$\phi_z(w) = (\phi_{z_1}(w_1), \phi_{z_2}(w_2))$$

in the bidisc $D^2$ for each point $z = (z_1, z_2) \in D^2$. For a fixed point $z \in D^2$, we define a unitary operator $U_z$ on $L^2(D^2)$ by

$$U_z h(w) = h \circ \phi_z(w) k_z(w)$$

for all $h \in H^2(D^2)$. It is easy to see that $U_z^* T_f U_z = T_{f \circ \phi_z}$. Thus

$$T_{f \circ \phi_z} T_{g \circ \phi_z} - T_{(f \circ \phi_z)(g \circ \phi_z)} = U_z^*(T_f T_g - T_{fg})U_z.$$

Therefore if $T_f T_g - T_{fg}$ is compact, then $T_{f \circ \phi_z} T_{g \circ \phi_z} - T_{(f \circ \phi_z)(g \circ \phi_z)}$ is compact.

Replacing f and g respectively by $f \circ \phi_{(z_1, 0, \cdots, 0)}$ and $g \circ \phi_{(z_1, 0, \cdots, 0)}$ in above analysis, we get that

$$\overline{\frac{\partial (f \circ \phi_{(z_1, 0)})}{\partial \bar{z}_1}(0, u_2)} \frac{\partial (g \circ \phi_{(z_1, 0)})}{\partial \bar{z}_1}(0, u_2) = 0$$

for all $z_1 \in D$ and $u_2 \in T$. But

$$\overline{\frac{\partial (f \circ \phi_{(z_1, 0, \cdots, 0)})}{\partial \bar{z}_1}(0, u_2)} = (|z_1|^2 - 1)\overline{\frac{\partial f}{\partial \bar{z}_1}(z_1, u_2)}.$$

Therefore

$$\overline{\frac{\partial f}{\partial \bar{z}_1}(z_1, u_2)} \frac{\partial g}{\partial \bar{z}_1}(z_1, u_2) = 0$$

for $u_2 \in T$ and $z_1 \in D$. This completes the proof of the theorem except we still need to prove our claim.

We first discuss $H_1(z)$. Let

$$(9) \qquad A(z) := \int_{T^2} \overline{f(U_\theta z)} g(U_\theta z) e^{i(m,\theta)} d\theta$$

for $z \in D^2$. Since both f and g are in $L^\infty(T^2)$, we write the harmonic extension of $f$ and $g$ in $D^2$ as

$$f = \sum_{l \in Z \times Z} f_l z^l \quad \text{and} \quad g = \sum_{j \in Z \times Z} g_j z^j$$

where $f_l$ and $g_j$ are two sequences of numbers satisfying

$$\sum_{l \in Z \times Z} |f_l|^2 < \infty \quad \text{and} \quad \sum_{j \in Z \times Z} |g_j|^2 < \infty.$$

A straightforward computation gives that

$$(10) \qquad A(z) = \sum_{j \in Z \times Z} \overline{f_{j+m}} g_j z^{-(j+m)} z^j$$

for $z \in D^2$, where recall that by our convention $z_1^{-1}$ is $\bar{z}_1$, for example. Since

$$\sum_{j \in Z \times Z} |\overline{f_{j+m}} g_j| \leq \left( \sum_{j \in Z \times Z} |f_{j+m}|^2 \right)^{1/2} \left( \sum_{l \in Z \times Z} |g_j|^2 \right)^{1/2} < \infty,$$

$A(z)$ is continuous on the closure $\overline{D}^2$ of $D^2$; and indeed $A(z)$ is defined by the series expansion in (10) for all $z \in \overline{D}^2$. But note that here $A(z)$ for $z \in D^2$ is not necessary the harmonic extension of $A(w)$ for $w \in T^2$. By an abuse of notation, we have that (9) holds for all $z \in \overline{D}^2$. It follows from the definition of $H_1(z)$ that for all $z \in D^2$,

$$
\begin{aligned}
H_1(z) &= -\int_{T^2} \int_{T^2} \overline{f(w)} g(w) |k_{U_{\theta}z}(w)|^2 d\sigma(w) e^{i(m,\theta)} d\theta \\
&= -\int_{T^2} \int_{T^2} \overline{f(U_\theta w)} g(U_\theta w) e^{i(m,\theta)} d\theta |k_z(w)|^2 d\sigma(w) \\
&= -\int_{T^2} A(w) |k_z(w)|^2 dA(w)
\end{aligned}
$$

where the second equality is obtained by changing the order of integration and a change of variables from $w$ to $U_\theta w$. That is, $H_1(z)$ in $D^2$ is the harmonic extension of the continuous function $A(w)$ on $T^2$. Therefore $H_1(z)$ is continuous on $\overline{D}^2$. Let $\phi_{z_2}(\lambda_2)$ be the Möbius map $\dfrac{z_2 - \lambda_2}{1 - \overline{z_2}\lambda_2}$. Changing variable $w_2 = \phi_{z_2}(\lambda_2)$ in the above integral, we have

$$H_1(z) = -\int_T \int_T A(w_1, \phi_{z_2}(\lambda_2))|k_{z_1}(w_1)|^2 dA(w_1)dA(\lambda_2).$$

It follows from Lebesque dominated convergence theorem and the continuity of $A(z)$ on $\overline{D^2}$ that

$$H_1(z, u_2) = \lim_{(z)_\alpha \to (z_1, u_2)} H_{22}(z) = -\int_T A(w_1, u_2)|k_{z_1}(w_1)|^2 dA(w_1).$$

That is, $H_1(z, u_2)$ as a function of $z_1$ on $D$ is the harmonic extension of $A(w_1, u_2)$ as a function of $w_1$ on $T$. Therefore $H_1(z, u_2)$ is harmonic in $z_1 \in D$ for any $u_2 \in T$.

Next we discuss $H_{12}(z)$. Without loss of generality, we consider a term $B(z)$ of $H_{12}$ as in (8).

$$\begin{aligned} B(z) &= \int_{T^2} \int_{T^2} \overline{f_{-+}(e^{i\theta_1}z_1, w_2)} g_{++}(w)|k_{U_\theta z}(w)|^2 dA(w)e^{i(m,\theta)}d\theta \\ &= \int_{T^2} \int_{T^2} \overline{f_{-+}(U_\theta(z_1, w_2))} g_{++}(U_\theta w)e^{i(m,\theta)}d\theta |k_z(w)|^2 dA(w) \end{aligned}$$

where the second equality is obtained by changing the order of integration and a change of variables from $w$ to $U_\theta w$. By integrating with respect to $w_1$ and noting that $f_{-+}(U_\theta(z_1, w_2))$ is independent of $w_1$, we get that

$$B(z) = \int_T \int_{T^2} \overline{f_{-+}(U_\theta(z_1, w_2))} g_{++}(U_\theta(z_1, w_2))e^{i(m,\theta)}d\theta |k_{z_2}(w)|^2 dA(w_2).$$

Now we define

$$(11) \qquad C(z) = \int_{T^2} \overline{f_{-+}(z_1, w_2)} g_{++}(w)e^{i(m,\theta)}d\theta.$$

By the proof of the continuity of $A(z)$, we see that $C(z)$ is continuous on $\overline{D^2}$. By an abuse of notation we see that (11) holds for all $z \in \overline{D^2}$. Indeed we have

$$(12) \qquad C(z) = \sum_{j \in Z \times Z} \overline{a_{j+m}} b_j z^{-(j+m)} z^j,$$

for all $z \in \overline{D^2}$, if we write $f_{-+}(z)$ and $g_{++}(z)$ as

$$f_{-+}(z) = \sum_{l \in Z \times Z} a_l z^l \quad \text{and} \quad g_{--}(z) = \sum_{j \in Z \times Z} b_j z^j$$

where, since $f_{-+}(z)$ and $g_{--}(z)$ are in $L^2(T^2)$, $a_l$ and $b_j$ are two sequences of numbers satisfying

$$\sum_{l \in Z \times Z} |a_l|^2 < \infty, \quad \text{and} \quad \sum_{j \in Z \times Z} |b_j|^2 < \infty.$$

Next note that

$$B(z) = \int_T C(z_1, w_2)|k_{z_2}(w)|^2 dA(w_2).$$

Hence $B(z)$ is continuous on $\overline{D^2}$. Let $\phi_{z_2}(\lambda_2)$ be the Möbius map $\dfrac{z_2 - \lambda_2}{1 - \overline{z_2}\lambda_2}$. We have

$$B(z) = \int_T C(z_1, \phi_{z_2}(\lambda_2)) dA(\lambda_2).$$

It follows from the continuity of $C(z)$ that

$$B(z_1, u_2) = \lim_{(z)_\alpha \to (z_1, u_2)} B(z) = \int_T C(z_1, u_2) dA(\lambda_2) = C(z_1, u_2).$$

To prove that $C(z_1, u_2)$ is harmonic in $z_1$ for any $u_2 \in T$, we look at the series expansion of $C(z)$ more carefully. We first note that $a_l = 0$ for $l = (l_1, l_2) \in Z \times Z$ such that $l_1 \geq 0$ and $b_j = 0$ except for $j = (j_1, j_2) \in Z_+ \times Z_+$. Therefore in fact $C(z_1, u_2) = 0$ for $m = (m_1, m_2) \in Z_+ \times Z$, and

$$C(z_1, u_2) = z_1^{-m_1} u_2^{m_2} \left( \sum_{0 \leq j_1 < -m_1} \sum_{j_2 \in Z_+} \overline{a_{j+m}} b_j \right)$$

for $m = (m_1, m_2) \in Z_- \times Z$. This finishes the proof of the claim concerning $H_{21}(z)$.

The proof of the claim concerning $H_{11}(z)$ and $H_{21}(z)$ is similar. So does the proof of the continuity of $H_{22}(z)$. Nevertheless one does not necessarily have that $H_{22}(z_1, u_2)$ is harmonic in $z_1$. So the claim is established. This completes the proof of the theorem.

We remark that in general the conditions (1) and (2) in Theorem 2 are not sufficient for $T_f T_g - T_{fg}$ to be compact. See Theorem 3 in next section for examples.

Next we use Theorem 2 to derive a result in [4] about compact Hankel operators on the bidisc; see [4] and reference therein for more related results.

**Corollary 1.** [4] *Let $f$ be a bounded functions on the bitorus $T^2$. The Hankel operators $H_f$ is compact if and only if it is zero.*

**Proof.** Note that $-H_f^* H_f = T_f T_f - T_{\bar f f}$. Thus by Theorem 2, $H_f$ is compact implies that for all $z_1 \in D$ and all $z_2 \in T$,

$$|\frac{\partial f}{\partial \bar z_1}(z_1, z_2)|^2 = 0$$

and for all $z_2 \in D$ and all $z_1 \in T$,

$$|\frac{\partial f}{\partial \bar z_2}(z_1, z_2)|^2 = 0.$$

Therefore $f$ is analytic in $D^2$. So $H_f$ is zero. The proof is complete.

The following result is an immediate consequence of Theorem 2. See also [11] for related results on the Bergman space of the bidisc.

**Corollary 2.** *Let $f$ and $g$ be two bounded functions on $T^2$. If one of the functions $f$, $\bar{f}$, $g$, $\bar{g}$ is analytic, then $T_f T_g - T_g T_f$ is zero if and only if $T_f T_g - T_g T_f$ is compact.*

**Proof.** Without loss of generality, we assume that $f$ is analytic. It is easy to see that in this case $T_f T_g - T_g T_f = T_f T_g - T_{fg}$.

By Theorem 2, if $T_f T_g - T_{fg}$ is compact, then for all $z_1 \in D$ and $z_2 \in T$,

$$\frac{\partial f}{\partial z_1}(z_1, z_2) \frac{\partial g}{\partial \bar{z}_1}(z_1, z_2) = 0,$$

and for all $z_2 \in D$ and all $z_1 \in T$,

$$\frac{\partial f}{\partial z_2}(z_1, z_2) \frac{\partial g}{\partial \bar{z}_2}(z_1, z_2) = 0.$$

Since $f$ is analytic on the unit disk, we conclude that for each $i$ ($i = 1, 2$) either $\bar{f}$ or $g$ is analytic in $z_i$. By Theorem 1, this implies that $T_f T_g - T_{fg} = 0$. This completes the proof of the corollary.

We remark that the above corollary is not valid without the assumption that one of the functions $f$, $\bar{f}$, $g$, $\bar{g}$ is analytic. An example will be shown in next section.

## 4. Compact semi-commutators

In this section we will characterize when $T_f T_g - T_{fg}$ is compact for Toeplitz operators $T_f$ and $T_g$ with symbols of the form $f(z_1, z_2) = f_1(z_1) f_2(z_2)$ and $g(z_1, z_2) = g_1(z_1) g_2(z_2)$. This will show that the condition in Theorem 2 is in general not a sufficient condition for $T_f T_g - T_{fg}$ to be compact. Furthermore it provides us with nontrivial compact semi-commutators.

**Theorem 3.** *Let $f_1$, $f_2$, $g_1$ and $g_2$ be nonzero bounded functions on the unit circle $T$. Let $f$ and $g$ be bounded functions on $T^2$ of the form $f(z_1, z_2) = f_1(z_1) f_2(z_2)$ and $g(z_1, z_2) = g_1(z_1) g_2(z_2)$. Then $T_f T_g - T_{fg}$ is a nonzero compact operator on the Hardy space of the bidisc if and only if the following two conditions hold.*
*(1) $f_1(z_1)g_1(z_1) = 0$ and $f_2(z_2)g_2(z_2) = 0$ on $T$, and*
*(2) $\lim_{z \to T^2} \|H_{\bar{f}} k_z\|_2 \|H_g k_z\|_2 = 0$.*

**Proof.** We first prove the necessity part of the theorem. Assume that $T_f T_g - T_{fg}$ is compact. By Theorem 2, we have that for all $z_1 \in D$ and all $z_2 \in T$,

$$(13) \qquad \frac{\partial f_1(z_1)}{\partial z_1} \frac{\partial g_1(z_1)}{\partial \bar{z}_1} f_2(z_2) g_2(z_2) = 0$$

and for all $z_2 \in D$ and all $z_1 \in T$,

$$(14) \qquad \frac{\partial f_2(z_2)}{\partial z_2} \frac{\partial g_2(z_2)}{\partial \bar{z}_2} f_1(z_1) g_1(z_1) = 0.$$

We claim that this implies that Condition (1) holds. If Condition (1) does not hold. Without loss of generality, assume that $f_1 g_1$ is not zero. Then equation (14) gives that

$$\frac{\partial f_2}{\partial z_2}\frac{\partial g_2}{\partial \bar{z}_2} = 0$$

for all $z_2 \in D$. That is either $\bar{f}_2$ or $g_2$ is analytic in $z_2$. Say, $g_2$ is analytic in $z_2$. Thus $g_2 f_2$ is not zero since otherwise either $f_2$ or $g_2$ will be zero. Now equation (13) implies that

$$\frac{\partial f_1}{\partial z_1}\frac{\partial g_1}{\partial \bar{z}_1} = 0.$$

That is either $\bar{f}_1$ or $g_1$ is analytic in $z_1$. Say, $\bar{f}_1$ is analytic in $z_1$. Hence $\bar{f}$ is analytic in $z_1$ and $g$ is analytic in $z_2$. By Theorem 1, we get that $T_f T_g - T_{fg}$ is zero. This is a contradiction. Therefore Condition (1) holds.

Next we show that Condition (2) holds. We first need some notations. Let $P_1$ be the projection defined by

$$P_1(h) = \sum_{m=(m_1,m_2) \in Z_+ \times Z} h_m z^m \quad \text{if } h = \sum_{m=(m_1,m_2) \in Z \times Z} h_i(z_2) z^m \in L^2(T^2),$$

and similarly $P_2$ the projection defined by

$$P_2(h) = \sum_{m=(m_1,m_2) \in Z \times Z_+} h_m z^m \quad \text{if } h = \sum_{m=(m_1,m_2) \in Z \times Z} h_i(z_2) z^m \in L^2(T^2),$$

It is easy to see that $P = P_1 P_2 = P_2 P_1$. Let $h_n(z_1) \in H^2(D)$ such that $h_n$ converges weakly to zero in $H^2(D)$. For any $h(z_2)$ in $H^2(D)$. Let $H_n(z_1, z_2)$ be $h(z_2)h_n(z_1)$. Note that $H_n$ converges weakly to zero in $H^2(D^2)$ as well. Note also by Condition (1),

$$(T_f T_g - T_{fg})H_n(z_1, z_2) = T_f T_g H_n(z_1, z_2)$$

$$= P_2(f_2 P_2 g_2 h)(z_2) T_{f_1} T_{g_1} h_n(z_1) = -P_2(f_2 P_2 g_2 h)(z_2) H^*_{\bar{f}_1} H_{g_1} h_n(z_1)$$

and $P_2(f_2)P_2(g_2 h)(z_2)$ is not zero for some $h \in H^2(D)$. Thus if $T_f T_g - T_{fg}$ is compact, then $\|H^*_{\bar{f}_1} H_{g_1} h_n\|$ converges to zero. Therefore $H^*_{\bar{f}_1} H_{g_1}$ is compact on the Hardy space of the unit disk. The compactness of $H^*_{\bar{f}_2} H_{g_2}$ can be established similarly. It follows from Theorem 2 in [13] that

$$\lim_{z_1 \to T} \|H_{\overline{f_1}} k_{z_1}\|_2 \|H_{g_1} k_{z_1}\|_2 = 0,$$

and

$$\lim_{z_2 \to T} \|H_{\overline{f_2}} k_{z_2}\|_2 \|H_{g_2} k_{z_2}\|_2 = 0.$$

So

$$\lim_{z_1 \text{ or } z_2 \to T} \|H_{\overline{f_1}} k_{z_1}\|_2 \|H_{g_1} k_{z_1}\|_2 \|H_{\overline{f_2}} k_{z_2}\|_2 \|H_{g_2} k_{z_2}\|_2 = 0.$$

On the other hand,

$$\|H_{\overline{f}}k_z\|_2\|H_g k_z\|_2 = \|H_{\overline{f_1}}k_{z_1}\|_2\|H_{g_1}k_{z_1}\|_2\|H_{\overline{f_2}}k_{z_2}\|_2\|H_{g_2}k_{z_2}\|_2.$$

Hence

$$\lim_{z\to T^2}\|H_{\overline{f}}k_z\|_2\|H_g k_z\|_2 = 0$$

as it is equivalent to

$$\lim_{z_1\to T}\|H_{\overline{f_1}}k_{z_1}\|_2\|H_{g_1}k_{z_1}\|_2 = 0,$$

and

$$\lim_{z_2\to T}\|H_{\overline{f_2}}k_{z_2}\|_2\|H_{g_2}k_{z_2}\|_2 = 0.$$

We now turn to the proof of the sufficient part of the theorem. By Condition (1), we have $T_f T_g - T_{fg} = T_f T_g$. Let $\phi$ and $\psi$ be in $H^2(D^2)$. We observe that

$$T_f T_g \phi = P_1 P_2(f_1 f_2 P_1 P_2(g_1 g_2 \phi))$$

$$= P_1 f_1 P_1 g_1 P_2(f_2 P_2(g_2 \phi))$$

$$= P_1 f_1 (I - P_1) g_1 P_2(f_2(I - P_2)(g_2 \phi))$$

$$= P_1 f_1 P_2(f_2(I - P_2)g_1(I - P_2)(g_2 \phi)),$$

where the third equality follows from Condition (1). Therefore

$$< T_f T_g \phi, \psi > = < P_1 f_1 P_2(f_2(I - P_1)g_1(I - P_2)g_2 \phi), \psi >$$

$$= < (I - P_1)g_1(I - P_2)g_2 \phi, \overline{f_1 f_2}\psi >$$

$$= < (I - P_1)g_1(I - P_2)g_2 \phi, (I - P_1)\overline{f_1}(I - P_2)\overline{f_2}\psi > .$$

It is easy to see that $(I - P_1)g_1(I - P_2)g_2\phi$ is co-analytic in $z_1$ and $z_2$. So

$$\bigtriangledown_1 (I - P_1)g_1(I - P_2)g_2 \phi = \frac{\partial}{\partial \overline{z_1}}(I - P_1)g_1(I - P_2)g_2 \phi.$$

Since

$$\int_{\partial D}(I - P_1)g_1(I - P_2)g_2 \phi(z_1, z_2)d\sigma(z_1) = 0$$

and

$$\int_{\partial D}\frac{\partial}{\partial \overline{z_1}}[(I - P_1)g_1(I - P_2)g_2 \phi(z_1, z_2)]d\sigma(z_2) = 0,$$

using the Littlewood-Paley formula, we have

$$< (I - P_1)g_1(I - P_2)g_2\phi, (I - P_1)\overline{f_1}(I - P_2)\overline{f_2}\psi >$$

$$= \int_D \int_D < \bigtriangledown_1 \bigtriangledown_2 (I - P_1)g_1(I - P_2)g_2\phi, \bigtriangledown_1 \bigtriangledown_2 (I - P_1)\overline{f_1}(I - P_2)\overline{f_2}\psi > \times$$

$$\log\frac{1}{|z_1|}\log\frac{1}{|z_2|}dA(z_1)dA(z_2)$$

$$= \int_D \int_D < \bigtriangledown_1(I - P_1)g_1 \bigtriangledown_2 (I - P_2)g_2\phi, \bigtriangledown_1(I - P_1)\overline{f_1} \bigtriangledown_2 (I - P_2)\overline{f_2}\psi > \times$$

$$\log\frac{1}{|z_1|}\log\frac{1}{|z_2|}dA(z_1)dA(z_2)$$

$$= I_1 + I_2 + I_3$$

where

$$I_1 = \int_{rD}\int_{rD} < \bigtriangledown_1 \bigtriangledown_2 (I - P_1)g_1(I - P_2)g_2\phi, \bigtriangledown_1 \bigtriangledown_2 (I - P_1)\overline{f_1}(I - P_2)\overline{f_2}\psi >$$

$$\times \log\frac{1}{|z_1|}\log\frac{1}{|z_2|}dA(z_1)dA(z_2),$$

$$I_2 = \int_{D/rD}\int_D < \bigtriangledown_1 \bigtriangledown_2 (I - P_1)g_1(I - P_2)g_2\phi, \bigtriangledown_1 \bigtriangledown_2 (I - P_1)\overline{f_1}(I - P_2)\overline{f_2}\psi >$$

$$\times \log\frac{1}{|z_1|}\log\frac{1}{|z_2|}dA(z_1)dA(z_2),$$

$$I_3 = \int_{rD}\int_{D/rD} < \bigtriangledown_1 \bigtriangledown_2 (I - P_1)g_1(I - P_2)g_2\phi, \bigtriangledown_1 \bigtriangledown_2 (I - P_1)\overline{f_1}(I - P_2)\overline{f_2}\psi >$$

$$\times \log\frac{1}{|z_1|}\log\frac{1}{|z_2|}dA(z_1)dA(z_2).$$

It is easy to see that there is a compact operator $K_r$ such that

$$I_1 = < K_r\phi, \psi > .$$

To deal with $I_2$ and $I_3$, we use the distribution function inequality established in [13]. As in the proof of Theorem 7 in [13], we obtain that

$$| \int_{D/rD} < \bigtriangledown_1(I - P_1)g_1 \bigtriangledown_2 (I - P_2)g_2\phi, \bigtriangledown_1(I - P_1)\overline{f_1} \bigtriangledown_2 (I - P_2)\overline{f_2}\psi > \log\frac{1}{|z_1|}dA(z_1)|$$

$$\leq C \sup_{|z_1|>r} \|H_{\overline{f_1}}k_{z_1}\|\|H_{g_1}k_{z_1}\|\| \bigtriangledown_2 (I - P_2)g_2\phi\|_{z_1}\| \bigtriangledown_2 (I - P_2)\overline{f_2}\psi\|_{z_1},$$

where for $h \in L^2(T^2)$ and $z_2$ fixed, $\|h\|_{z_1}$ denotes the $L^2$ norm of $h(z_1, z_2)$ as a function of $z_1$. Similarly we also have

$$| \int_{D/rD} < \bigtriangledown_1 (I - P_1)g_1 \bigtriangledown_2 (I - P_2)g_2\phi, \bigtriangledown_1(I - P_1)\overline{f_1} \bigtriangledown_2 (I - P_2)\overline{f_2}\psi > \log \frac{1}{|z_2|} dA(z_2)|$$

$$\leq C \sup_{|z_2|>r} \| \|H_{\overline{f_2}}k_{z_2}\| \|H_{g_2}k_{z_2}\| \| \bigtriangledown_1 (I - P_1)g_1\phi\|_{z_2} \| \bigtriangledown_1 (I - P_1)\overline{f_1}\psi\|_{z_2}.$$

So

$$|I_2| \quad \leq C \sup_{|z_1|>r} \| H_{\overline{f_1}}k_{z_1}\| \|H_{g_1}k_{z_1}\| \times$$

$$\int_D \| \bigtriangledown_2 (I - P_2)g_2\phi\|_{z_1} \| \bigtriangledown_2 (I - P_2)\overline{f_2}\psi\|_{z_1} \log \frac{1}{|z_2|} dA(z_2).$$

By the Holder inequality,

$$|I_2| \quad \leq C \sup_{|z_1|>r} \|H_{\overline{f_1}}k_{z_1}\| \|H_{g_1}k_{z_1}\| \times$$

$$[\int_D \int_{\partial D} | \bigtriangledown_2 (I - P_2)g_2\phi(z_1, z_2)|^2 \log \frac{1}{|z_2|} d\sigma(z_1) dA(z_2)]^{1/2} \times$$

$$[\int_D \int_{\partial D} | \bigtriangledown_2 (I - P_2)\overline{f_2}\psi(z_1, z_2)|^2 \log \frac{1}{|z_2|} d\sigma(z_1) dA(z_2)]^{1/2}.$$

By the Littlewood-Paley theorem and the fact that $g_2$ and $f_2$ are bounded,

$$|I_2| \leq C \sup_{|z_1|>r} \|H_{\overline{f_1}}k_{z_1}\| \|H_{g_1}k_{z_1}\| \|\phi\| \|\psi\|.$$

Similarly we can also get the estimate for $I_3$,

$$|I_3| \leq C \sup_{|z_2|>r} \|H_{\overline{f_2}}k_{z_2}\| \|H_{g_2}k_{z_2}\| \|\phi\| \|\psi\|.$$

In the proof of the necessary part of the theorem we have shown that the second condition is equivalent to

$$\lim_{z_1 \to T} \|H_{\overline{f_1}}k_{z_1}\|_2 \|H_{g_1}k_{z_1}\|_2 = 0,$$

and

$$\lim_{z_2 \to T} \|H_{\overline{f_2}k_{z_2}}\|_2 \|H_{g_2}k_{z_2}\|_2 = 0.$$

Hence we conclude that

$$\lim_{r \to 1} \|T_f T_g - T_{fg} - K_r\| = 0.$$

Therefore $T_f T_g - T_{fg}$ is compact. This completes the proof of the theorem.

We now give an example to show that the assumption in Corollary 2 is necessary. Namely, there exist bounded functions $f$ and $g$ such that $T_f T_g - T_g T_f$ is compact, but $T_f T_g - T_g T_f$ is not zero.

**Example.** Let $f_1(z_1)$ and $g_1(z_1)$ be nonzero continuous functions on the unit circle satisfying the conditions that $f_1(z_1)g_1(z_1) = 0$ and the set $\{z_1 : f_1(z_1) = 0\} \cap \{z_1 : g_1(z_1) = 0\}$ as a subset of $T$ has positive Lebesgue measure. Let $f_2(z_2)$ and $g_2(z_2)$ be nonzero continuous functions on the unit circle such that $f_2(z_2)g_2(z_2) = 0$.

First we show that $T_f T_g - T_g T_f$ is not zero. Note that

$$(T_f T_g - T_g T_f)h_1(z_1) = P_2(f_2)P_2(g_2)(T_{f_1}T_{g_1} - T_{g_1}T_{f_1})h_1(z_1)$$

for all $h_1 \in H^2(D)$. It is easy to see that the assumption on $f_2$ and $g_2$ implies that $P_2(f_2)P_2(g_2)$ is not zero. Therefore it suffices to show that $T_{f_1}T_{g_1} - T_{g_1}T_{f_1}$ is not zero. Let $S_1$ denote the multiplication by $z_1$ on $H^2(D)$. We observe that $I - S_1 S_1^*$ is the rank one operator $e_0 \otimes e_0$, where $e_0 = 1 \in H^2(D)$. Note that

$$
\begin{aligned}
& T_{f_1}T_{g_1} - T_{g_1}T_{f_1} - S_1^*(T_{f_1}T_{g_1} - T_{g_1}T_{f_1})S_1 \\
= \; & T_{f_1}T_{g_1} - S_1^* T_{f_1}T_{g_1}S_1 - (T_{g_1}T_{f_1} - S_1^* T_{g_1}T_{f_1}S_1) \\
= \; & S_1^* T_{f_1}S_1 S_1^* T_{g_1}S_1 - S_1^* T_{f_1}T_{g_1}S_1 - (S_1^* T_{g_1}S_1 S_1^* T_{f_1}S_1 - S_1^* T_{g_1}T_{f_1}S_1) \\
= \; & -S_1^* T_{f_1}(I - S_1 S_1^*)T_{g_1}S_1 + S_1^* T_{g_1}(I - S_1 S_1^*)T_{f_1}S_1 \\
= \; & -S_1^* T_{f_1}e_0 \otimes S_1^* T_{g_1}^* e_0 + S_1^* T_{g_1}e_0 \otimes S_1^* T_{f_1}^*,
\end{aligned}
$$

where the second equality follows from the fact that for Toeplitz operator $T_{f_1}$, $S_1^* T_{f_1}S_1 = T_{f_1}$. Thus if $T_{f_1}T_{g_1} - T_{g_1}T_{f_1} = 0$, then there exists a contant $c$ such that

$$S_1^* T_{f_1}e_0 = cS_1^* T_{g_1}e_0.$$

That is, there exists a constant $d$ such that $P_1(f_1 - cg_1) = d$. Hence if we denote the function $\overline{f_1} - \overline{cg_1}$ by $h$, then $h$ is analytic. But by our assumption, $h(z_1) = 0$ on a subset of $T$ with positive Lebesgue measure. Therefore $h_1 = 0$ on $T$, i.e., $f_1 = cg_1$. This contradicts to the conditions that both $f_1$ and $g_1$ are nonzero and $f_1 g_1 = 0$ on $T$. So $T_{f_1}T_{g_1} - T_{g_1}T_{f_1}$ is not zero.

Next we show that $T_f T_g - T_g T_f$ is compact. Since $f_1(z_1)$, $g_1(z_1)$, $f_2(z_2)$ and $g_2(z_2)$ are continuous, so $H_{\overline{f_1}}H_{g_1}$ and $H_{\overline{f_2}}H_{g_2}$ are compact. Hence

$$\lim_{z_1 \to T} \|H_{\overline{f_1}}k_{z_1}\|_2 \|H_{g_1}k_{z_1}\|_2 = 0,$$

and

$$\lim_{z_2 \to T} \|H_{\overline{f_2}k_{z_2}}\|_2 \|H_{g_2}k_{z_2}\|_2 = 0,$$

or equivalently Condition (2) in Theorem 3 holds. Therefore $T_f T_g - T_{fg}$ and $T_g T_f - T_{gf}$ are compact. Note that

$$T_f T_g - T_g T_f = (T_f T_g - T_{fg}) - (T_g T_f - T_{gf}).$$

Hence $T_f T_g - T_g T_f$ is compact.

## References

[1] S. Axler, S.-Y. A. Chang, D. Sarason, *Product of Toeplitz Operators*, Integral Equations and Operator Theory, 1 (1978), 285-309.

[2] A. Brown and P. R. Halmos *Algebraic properties of Toeplitz operators*, J. Reine Angew. Math., 213 (1963), 89-102.

[3] S.-Y. A. Chang and R. Fefferman, *Some recent developments in Fourier analysis and $H^p$-theory on the product domain*, Bulletin of Amer. Math. Soc., 12 (1985), 1-43.

[4] M. Cotlar and C. Sadosky, *Abstract, Weighted, and multidimensional Adamjan-Arov-Krein theorem, and the singular numbers of Sarason commutants*, Integral Equations and Operator Theory, 17 (1993) 169-201.

[5] M. Cotlar and C. Sadosky, *Two distinguished subspaces of product BMO and the Nehari-AAK theory for Hankel operators on the torus*, preprint, MSRI, 1995.

[6] R. G. Douglas, *Banach algebra techniques in the operator theory*, Academic Press, New York and London, 1972.

[7] R. G. Douglas, *Banach algebra techniques in the theory of Toeplitz operators*, Regional Conference Series in Mathematics, American Mathematical Society 15, 1972.

[8] N. K. Nikolskii, Treatise on the shift operator, Springer-Verlag, New York, 1985.

[9] W. Rudin, Function theory on the polydisk, Benjamin Inc., New York, 1969.

[10] E. M. Stein and G. Weiss, Introduction to Fourier analysis on euclidean spaces, Princeton University Press 1971.

[11] S. Sun and D. Zheng, *Toeplitz operators on the polydisk*, to appear in Proc. Amer. Math. Soc..

[12] A. Volberg, *Two remarks concerning the theorem of S. Axler, S.-Y. A. Chang, and D. Sarason*, J. Operator Theory, 8 (1982) 209-218.

[13] D. Zheng, *The Distribution function inequality and Products of Toeplitz operators and Hankel operators*, to appear in J. Functional Analysis.

CAIXING GU, MATHEMATICS DEPARTMENT, UNIVERSITY OF CALIFORNIA IRVINE, IRVINE, CA 72717
    *E-mail address*: `cgu@math.uci.edu`

DECHAO ZHENG, MATHEMATICS DEPARTMENT, MICHIGAN STATE UNIVERSITY, EAST LANSING, MI 48824
    *E-mail address*: `dechao@math.msu.edu`